\documentclass[11pt]{article}
\usepackage{amsmath,mathrsfs,latexsym,amsfonts,theorem}

{\theorembodyfont{\slshape}
\newtheorem{theorem}{Theorem}
\newtheorem{proposition}{Proposition}
\newtheorem{lemma}{Lemma}
\newtheorem{corollary}{Corollary}
}
{\theorembodyfont{\rmfamily}
\newtheorem{definition}{Definition}

\newtheorem{remark}{Remark}

}

\def\proof{{\noindent\sc Proof. \quad}}
\newcommand\proofof[1]{{\medskip\noindent\sc Proof of #1.\quad}}
\def\eproof{\hfill\qed\medskip}
\newcommand\qed{{\unskip\nobreak\hfil\penalty50\hskip2em\vadjust{}
\nobreak\hfil$\Box$\parfillskip=0pt\finalhyphendemerits=0\par}}

\renewcommand{\hat}{\widehat}

\makeatletter
\def\R{\mathchoice{{\setbox0=\hbox{\rm I}\copy0\kern-0.55\wd0\hbox{\rm R}}}{%
  {\setbox0=\hbox{\rm I}\copy0\kern-0.55\wd0\hbox{\rm R}}}{%
  {\setbox0=\hbox{$\m@th\scriptstyle\rm I$}\copy0\kern-0.5\wd0%
  \scriptstyle\rm R}}{%
  {\setbox0=\hbox{$\m@th\scriptscriptstyle\rm I$}\copy0\kern-0.6\wd0%
  \scriptscriptstyle\rm R}}}
\def\N{\mathchoice{{\setbox0=\hbox{\rm I}\copy0\kern-0.55\wd0\rm N}}{%
  {\setbox0=\hbox{\rm I}\copy0\kern-0.55\wd0\rm N}}{%
  {\setbox0=\hbox{$\m@th\scriptstyle\rm I$}\copy0\kern-0.5\wd0%
  \scriptstyle\rm N}}{%
  {\setbox0=\hbox{$\m@th\scriptscriptstyle\rm I$}\copy0\kern-0.6\wd0%
  \scriptscriptstyle\rm N}}}
\def\Z{\mathchoice{{\setbox0=\hbox{\rm Z}\copy0\kern-0.5\wd0\box0}}{%
  {\setbox0=\hbox{\rm Z}\copy0\kern-0.5\wd0\box0}}{%
  {\setbox0=\hbox{$\m@th\scriptstyle\rm Z$}\copy0\kern-0.55\wd0\box0}}{%
  {\setbox0=\hbox{$\m@th\scriptscriptstyle\rm Z$}\copy0\kern-0.6\wd0\box0}}}
\def\Q{\mathchoice{{\setbox0=\hbox{\rm Q}\kern0.2\wd0 0\kern-0.9\wd0\box0}}{%
  {\setbox0=\hbox{\rm Q}\kern0.2\wd0 0\kern-0.9\wd0\box0}}{%
  {\setbox0=\hbox{$\m@th\scriptstyle\rm Q$}\kern0.2\wd0 0\kern-0.9\wd0\box0}}{%
  {\setbox0=\hbox{$\m@th\scriptscriptstyle\rm Q$}\kern0.2\wd0 0%
  \kern-0.9\wd0\box0}}}
\def\C{\mathchoice{{\setbox0=\hbox{\rm C}\kern0.47\wd0%
  {\vrule height0.95\ht0depth1.1\dp0width0.7pt}\kern-0.47\wd0\box0}}{%
  {\setbox0=\hbox{\rm C}\kern0.47\wd0%
  {\vrule height0.95\ht0depth1.1\dp0width0.7pt}\kern-0.47\wd0\box0}}{%
  {\setbox0=\hbox{$\m@th\scriptstyle\rm C$}\kern0.47\wd0%
  {\vrule height0.95\ht0depth1.1\dp0width0.5pt}\kern-0.47\wd0\box0}}{%
  {\setbox0=\hbox{$\m@th\scriptscriptstyle\rm C$}\kern0.47\wd0%
  {\vrule height0.95\ht0depth1.1\dp0width0.4pt}\kern-0.47\wd0\box0}}}

\def\Prob{\mathop{\mathsf{Prob}}}

\def\Oh{{\cal O}}

\def\a{{\alpha}}

\def\g{{\gamma}}
\def\BB {{\cal B}}

\def\DD {{\cal D}}
\def\Q {{\cal Q}\,}

\def\Z {{\cal Z}\,}
\def\Id{{\mathrm{Id}}}

\def\m{{\boldsymbol m}}
\def\c{{\boldsymbol c}}
\def\ckl{{\boldsymbol c}^{\dagger}_{k\ell}}
\def\ck{{\boldsymbol c}_{k}}
\def\bfE{{\mathbf E}}

\def\cd{{\boldsymbol c}_{\mathsf{det}}}
\def\MM{{\mathscr M}}
\def\RR{{\mathscr R}}

\def\rt{^{\mathsf T}}

\def\JACM{Journal of the ACM}

\begin{document}

\begin{title}
 {\LARGE \bf Componentwise condition numbers of random sparse matrices}
\end{title}
\author{Dennis Cheung\\
United International College\\
Tang Jia Wan\\
Zhuhai, Guandong Province\\
P.R. of CHINA\\
e-mail:
{\tt dennisc@uic.edu.hk}
\and
Felipe Cucker
\thanks{Partially supported by 
GRF grant CityU 100808}
\\
Department of Mathematics\\
City University of Hong Kong\\
83 Tat Chee Avenue, Kowloon\\
HONG KONG\\
e-mail:
{\tt macucker@cityu.edu.hk}
}

\date{}

\maketitle

\begin{quote}
{\footnotesize {\bf Abstract.\quad}
We prove an $\Oh(\log n)$ bound for the expected value of the logarithm of 
the componentwise (and, a fortiori, the mixed) condition number of a random 
sparse $n\times n$ matrix. As a consequence, small bounds on the average 
loss of accuracy for triangular linear systems follow.  
}
\end{quote}

\section{Introduction}

Triangular systems of linear equations provide one of the few examples 
in numerical linear algebra where a gap occurs between stability analysis 
and everyday practice. One could summarize this gap as follows:
\begin{quote}
{\sl Triangular systems of equations are generally solved to high accuracy 
in spite of being, in general, ill-conditioned.} 
\end{quote}
This state of affairs had been already noted by J.H.~Wilkinson 
in~\cite[p.~105]{Wilkinson63}: ``In practice one almost invariably finds 
that if $L$ is ill-conditioned, so that $\|L\|\|L^{-1}\|\gg1$, then the 
computed solution of $Lx=b$ (or the computed inverse) is far more accurate 
than [what forward stability analysis] would suggest.'' 

An explanation to this gap is suggested by N.~Higham~\cite{Higham:89} who 
notes that the backward error analysis given by Wilkinson for the solution 
of triangular systems yields (small) {\em  componentwise} bounds on the perturbation 
matrix (see Section~\ref{sec:final}(1) below). 
Higham then uses this fact to deduce small forward error bounds 
for particular subclasses of triangular systems and to numerically investigate 
the accuracy of other particular such systems. In doing so, Higham makes 
use of the mixed condition number introduced by 
Skeel~\cite{skeel}\footnote{Skeel called it ``componentwise.'' In this 
paper, however, following the notation introduced in~\cite{Gohberg}, we will 
use this word for the condition numbers measuring both data perturbation 
and computed errors in a componentwise fashion.}. This condition number 
has a natural role in analyzing accuracy of triangular systems since bounds 
for it, together with the backward analysis of Wilkinson mentioned above, 
yield forward analysis bounds for the computed solution of the system.  
Furthermore, the restriction to componentwise perturbations ---both 
in the backward error analysis and in the mixed condition number--- forces 
perturbations to preserve the triangular structure of the data matrices. 

A further step in explaining the gap, somehow orthogonal to the work 
of Higham, was given by D.~Viswanath and N.~Trefethen in~\cite{ViTr98}
where a precise meaning to the expression ``triangular systems are, 
in general, ill-conditioned'' was given. Indeed, if $L_n$ denotes a 
random triangular $n\times n$ matrix (whose entries are independent 
standard Gaussian random variables) and 
$\kappa_n=\|L_n\|\|L_n^{-1}\|$ is its condition number (which is a 
positive random variable) then, it is shown in~\cite{ViTr98},
$$
   \sqrt[n]{\kappa_n}\to 2 \quad\mbox{almost surely}
$$ 
as $n\to\infty$. A straightforward consequence of this result is that 
the expected value of $\log\kappa_n$ satisfies 
$\bfE(\log\kappa_n)=\Omega(n)$. 

The goal of this paper is to close the gap by giving a precise meaning, 
in the sense of~\cite{ViTr98}, to the other half of the statement above 
namely, to the expression ``triangular systems are generally solved to 
high accuracy.'' More precisely, we consider the mixed condition numbers 
$\m^{\dagger}(L_n)$ and $\m(L_n,b_n)$ for the problems of matrix 
inversion and linear equation solving, respectively, for a random 
triangular $L_n$ as above and a random $b_n\in\R^n$. Then, we show that 
\begin{equation}\label{eq:main}
   \bfE(\log \m^{\dagger}(L_n)),\,
   \bfE(\log \m(L_n,b_n)) =\Oh(\log n). 
\end{equation}
From the bound on $\bfE(\log \m(L_n,b_n))$ it follows that the average loss 
of precision in the solution of random triangular systems is small. From 
that on $\bfE(\log \m^{\dagger}(L_n))$, that the one for matrix inversion is 
small as well. One can therefore replace the summary above by the following:
\begin{quote}
{\sl Triangular systems of equations are generally solved to high accuracy 
because their backward error analysis yields small componentwise 
perturbations and triangular matrices are, in general, well conditioned 
for these perturbations.  
} 
\end{quote}

The results showing~\eqref{eq:main}, Theorems~\ref{th3} and~\ref{th5} 
below, are proved in the more general context of sparse matrices 
(i.e.~matrices with a fixed pattern of zeros) and componentwise 
condition numbers (which ensure high relative accuracy {\em in each component} 
of the computed solution $A^{-1}$ or $x$). Besides triangular 
matrices, these results apply to other classes of sparse matrices such
as, for instance, tridiagonal matrices. In the process of proving them, 
we found useful to estimate as well the average mixed condition for the 
computation of the determinant.  
\bigskip

\section{Preliminaries} 

Condition numbers measure the worst-case magnification of a small 
data perturbation 
in the computed outcome. As originally introduced by Turing~\cite{Turing48}, 
they were {\em normwise} in the sense that both the data perturbation and the 
outcome's error are measured using norms (in the space of data and outcomes 
respectively). In contrast, {\em mixed} condition numbers measure data 
perturbation componentwise, and {\em componentwise} condition numbers measure 
both data perturbation and outcome's error in this way. 

To define these condition numbers the following
form of ``distance'' function will be useful.
For points $u,v \in \R^p$ we define 
$\frac{u}{v}=(w_1,\ldots,w_p)$ with 
$$
w_i=\left\{ \begin{array}{ll}
             u_i/v_i &\mbox{if $v_i\neq 0$}\\ 
             0 &\mbox{if $u_i=v_i=0$}\\ 
             \infty &\mbox{otherwise.}
             \end{array}\right.
$$
Then we define
$$
   d(u,v)=\left\|\frac{u-v}{v}\right\|_\infty.
$$
Note that, if $d(u,v)<\infty$, 
$$
  d(u,v)=\min\{\nu \geq 0\mid |u_i-v_i|\leq\nu|v_i|
              \mbox{ for $i=1,\ldots,p$}\}.
$$
For $\delta>0$ and $a\in\R^p$ we denote
$\BB(a,\delta)=\{x\in\R^p\mid d(x,a)\leq \delta\}$. 

\begin{definition}\label{def:CN}
Let $\DD\subseteq\R^p$ and $F:\DD \rightarrow \R^q$ be a continuous
mapping. Let $a \in \DD$ such that $F(a)\neq 0$.
\begin{description}
\item[(i)] The {\em mixed condition number} of $F$ at $a$  (with respect to 
a norm $\|\ \|_q$ on  $\R^q$) is
defined by
$$
  \m(F,a)=\lim_{\delta \rightarrow 0} \sup_{x \in \BB(a,\delta) \atop x \neq
  a}\frac{\|F(x)-F(a)\|_q}{\|F(a)\|_q}\frac{1}{d(x,a)}.
$$
\item[(ii)]
Suppose $F(a)=(f_1(a),\ldots,f_q(a))$ is such
that $f_j(a)\neq 0$ for $j=1,\ldots,q$. Then the {\em componentwise
condition number} of $F$ at $a$ is
$$
  \c(F,a)=\lim_{\delta \rightarrow 0} 
    \sup_{x \in \BB(a,\delta) \atop x \neq a}
    \frac{d(F(x),F(a))}{d(x,a)}.
$$
\end{description}
\end{definition}

\begin{proposition}\label{prop:monot}
For all $a\in\DD$ and any monotonic norm in $\R^q$, 
$\m(F,a)\leq \c(F,a)$.
\end{proposition}

\proof
For all $x\in\BB(a,\delta)$ and all $i\leq q$, 
$|F(x)_i-F(a)_i|\leq d(F(x),F(a))|F(a)_i|$. 
Since $\|\ \|$ is monotonic (cf.~\cite{BSW:61}), 
this implies $\|F(x)-F(a)\|\leq d(F(x),F(a))\|F(a)\|$ 
and hence the statement.
\eproof

In all what follows, for $n\in\N$, 
we denote the set $\{1,\ldots, n\}$ by $[n]$.

\begin{definition}
We denote by $\MM$ the set of $n\times n$ real matrices and by 
$\Sigma$ its subset of singular matrices. 
Also, for a subset $S\subseteq[n]^2$ 
we denote 
$$
  \MM_S=\{A\in\MM \mid \mbox{ if 
  $(i,j)\not\in S$ then $a_{ij}=0$}\}. 
$$ 
We denote by $\RR_S$ the space of random $n\times n$ matrices obtained by 
setting $a_{ij}=0$ if $(i,j)\in S$ and drawing all other entries 
independently from the standard Gaussian $N(0,1)$. As above, if 
$S=[n]^2$, we write simply $\RR$.  
\end{definition}

In the rest of this paper, for non-singular matrices $A,A'$, we denote 
their inverses by $\Gamma,\Gamma'$, respectively. Also, we denote 
by $A_{(ij)}$ the sub-matrix of $A$
obtained by removing from $A$ its $i$th row and its $j$th column.

\section{Determinant Computation}

We consider here the problem of computing the determinant of a matrix $A$ 
and its componentwise condition number $\cd(A)$. The main result of 
this section is the following. 

\begin{theorem}\label{s2}
For $S\subseteq[n]^2$  and $t\geq 2|S|$ we have 
$$
\Prob\{\cd(A)\geq t\}
 \leq |S|^2\frac{1}{t}.  
$$
\end{theorem}

Average loss of precision (in a base $b$) is measured by the expected 
value of the logarithm (in that base) of the condition number. 
We may use Theorem~\ref{s2} to obtain one such result for the computation 
of the determinant. To avoid problems caused by this condition number 
being less than 1 we consider the function $\log_+(x)$ defined to be 
$\log(x)$ if $x\geq 1$ and $0$ otherwise. 

\begin{corollary}\label{cor:s2}
For a base $b>1$, 
$\bfE(\log_+ \cd(A))\leq 
 2\log |S| +\frac{1}{\ln b}$ where 
$\bfE$ denotes expectation over $A\in\RR_S$.  
\end{corollary}

Towards the proof of Theorem~\ref{s2} we first 
obtain explicit expressions for $\cd(A)$. We begin by 
noting that taking 
$F:\MM\to \R$ to be $F(A)=\det(A)$ in Definition~\ref{def:CN} we 
obtain, for $A\in\MM\setminus\Sigma$, 
$$
  \cd(A)= \lim_{\delta\rightarrow 0} \sup_{A'\in \BB(A, \delta)}
         \frac{|\det(A')-\det(A)|}{\delta|\det(A)|}.
$$
Also, for $A\in\Sigma$, we have $\cd(A)=0$ if
$$
  \lim_{\delta\rightarrow 0} \sup_{A'\in \BB(A, \delta)}
           \frac{|\det(A')|}{\delta}=0
$$
and $\cd(A)=\infty$ otherwise. Note that $\cd(0)=0$. 

\begin{lemma}\label{s3}
For $A\in\MM\setminus\Sigma$,  
$$
   \cd(A)= \sum_{i,j\in[n]}|a_{ij}\g_{ji}|.
$$
For $A\in\Sigma$, $\cd(A)=0$ if
$$
  \sum_{i,j\in[n]}|a_{ij}\,\det (A_{(ij)})|=0
$$
and $\cd(A)=\infty$ otherwise.
\end{lemma}

\proof
Let $A\in\MM$. For any $i\in[n]$, expanding by the $i$th row, 
$$
 \det(A) = \sum_{j\in[n]} (-1)^{i+j} a_{ij}\, \det(A_{(ij)}).
$$
Hence, for all $i,j\in[n]$, 
$$
 \frac{\partial \det(A)}{\partial a_{ij}} =\det(A_{(ij)}).
$$
Using Taylor's expansion and these equalities we obtain
$$
 \det(A') = \det(A) + \sum_{i,j\in[n]}
  (a'_{ij}-a_{ij})\det(A_{(ij)}) + \Oh\left(\|A'-A\|^2\right).
$$
Here, the norm in $\|A'-A\|$ is not relevant since all 
norms in $\MM$ are equivalent. By choosing a monotonic norm we 
have that if $A'\in\BB(A,\delta)$ then $\|A'-A\|=\Oh(\delta)$. 
It follows that, for $A\not\in\Sigma$,  
\begin{eqnarray*}
  \cd(A) &=&\lim_{\delta\rightarrow 0}
            \sup_{A'\in\BB(A,\delta)}
            \frac{|\det(A')-\det(A)|}{\delta|\det(A)|}\\
       &=&\lim_{\delta\rightarrow 0} \sup_{A'\in\BB(A,\delta)}
           \sum_{i,j\in[n]}\frac{|(a'_{ij}-a_{ij})\det(A_{(ij)})|}
                {\delta|\det(A)|}\\
       &=&\lim_{\delta\rightarrow 0} \sup\left(
               \sum_{i,j\in[n]}\frac{|(a'_{ij}-a_{ij})\det(A_{(ij)})|}
               {\delta|\det(A)|} : \frac{|a'_{ij}-a_{ij}|}{|a_{ij}|}\leq\delta
               \right).
\end{eqnarray*}
The supremum above is attained by taking $a'_{ij}=a_{ij}(1+\delta)$ and 
therefore
\begin{eqnarray*}
  \cd(A) &=& \sum_{i, j\in[n]}
            \left|\frac{a_{ij}\det(A_{ij})}{\det(A)}\right|
         =\sum_{i, j\in[n]}|a_{ij}\g_{ji}|.
\end{eqnarray*}
If $A\in\Sigma$ it similarly follows that 
\begin{equation*}
 \lim_{\delta\rightarrow 0}
 \sup_{A'\in\BB(a,\delta)}\frac{|\det(A')|}{\delta}=\sum_{i,j\in[n]}
       |a_{ij}\,\det(A_{(ij)})|
\end{equation*}
and hence the statement.
\eproof

\begin{lemma}\label{sTail}
Let $p, q$ be two fixed vectors in $\R^n$ such that
$\|p\|\leq\|q\|$. If $x \sim N(0,\Id_n)$ then, for all $t\geq 2$, 
$$
   \Prob\left\{\left|\frac{x\rt p}{x\rt q}\right|\geq t\right\}
   \leq \frac{1}{t}.
$$
\end{lemma}

\proof 
Let $\nu=\|q\|$. By the orthogonal invariance of $N(0,\Id_n)$ we may assume 
$q=(\nu,0,\ldots,0)$. Also, by appropriately scaling, we may assume that 
$\nu=1$. Note that then, $\|p\|\leq 1$. We therefore have
\begin{eqnarray}\label{eq:Tail}
\Prob\left\{\left|\frac{x\rt p}{x\rt q}\right|\geq t \right\}
&=& \Prob\left\{\left|p_1+
 \sum_{i\in\{2,\ldots, n\}} 
 \frac{x_ip_i}{x_1}\right|\geq t\right\}\notag\\
&=& \Prob\left\{\left|p_1+\frac{1}{x_1} \a Z\right|\geq t\right\}\\
&=&  
  \Prob\left\{ \frac{Z}{x_1} \geq \frac{t-p_1}{\a}\right\} 
+   \Prob\left\{ \frac{Z}{x_1} \leq \frac{-t-p_1}{\a}\right\} \notag
\end{eqnarray}
where $Z=N(0,1)$ independent of $x_1$ and 
$\a=\sqrt{p_2^2+\ldots+p_n^2}\leq 1$. Here we used that a sum of 
independent centered Gaussians is a centered Gaussian whose variance 
is the sum of the terms variances. Note that in case $\a=0$ the statement 
is trivially true.

Since the $x_1$ and $Z$ are independent $N(0,1)$,  
the angle $\theta=\arctan\left(Z/x_1\right)$ 
is uniformly distributed in $[-\pi/2, \pi/2]$ and we have, for 
$\gamma\in\R$, 
\begin{eqnarray*}
    \Prob\left\{\frac{Z}{x_1} \geq \gamma \right\}
 &=& \Prob \left\{\theta \geq \arctan\gamma\right\}
     =\frac{1}{\pi}\left(\frac{\pi}{2}-\arctan \gamma\right)\\
 &=& \frac{1}{\pi}\int_{\gamma}^\infty\frac{1}{1+t^2} dt 
     \;\leq\; \frac{1}{\pi} \int_{\gamma}^\infty\frac{1}{t^2} dt
     \;=\; \frac{1}{\pi\gamma}.
\end{eqnarray*}
Similarly, for $\sigma\in\R$, 
\begin{eqnarray*}
    \Prob\left\{\frac{Z}{x_1} \leq \sigma \right\}
 &=& 1 - \Prob \left\{\theta \geq \arctan \sigma\right\}
     = 1 - \frac{1}{\pi}\left(\frac{\pi}{2}-\arctan \sigma\right)\\
 &=& = \frac{1}{\pi}\left(\frac{\pi}{2}-\arctan (-\sigma)\right)
      \;\leq\; \frac{1}{\pi(-\sigma)}.
\end{eqnarray*}
Using these bounds in~\eqref{eq:Tail} with $\gamma=\frac{t-p_1}{\a}$ 
and $\sigma=\frac{-t-p_1}{\a}$ we obtain
$$
 \Prob\left\{\left|\frac{x\rt p}{x\rt q}\right|\geq t \right\}
 \leq \frac{1}{\pi}\left(\frac{\a}{t-p_1}+\frac{\a}{t+p_1}\right)
 =\frac{\a}{\pi}\frac{2t}{t^2-p_1^2}
  \leq \frac{2}{\pi}\frac{t}{t^2-1}\leq \frac{1}{t}
$$
the last since $t\geq 2$. 
\eproof

\begin{lemma}\label{s23}
Let $S\subseteq[n]^2$ be such that $\MM_S\subseteq \Sigma$. 
Then, for all $A\in\MM_S$, $\cd(A) = 0$.
\end{lemma}

\proof 
Since $\MM_S\subseteq\Sigma$ and $A\in\MM_S$ we have 
$\BB(A,\delta)\subseteq\Sigma$ for all $\delta>0$. The result 
now follows. 
\eproof

\begin{lemma}\label{s22}
Let $S\subset[n]^2$ such that $\MM_S\not\subseteq\Sigma$. Then 
$$
  \Prob_{A\in\RR_S}(A\mbox{ is singular})= 0.
$$
\end{lemma}

\proof 
The set of singular matrices in $\MM_S$ is the zero set of the 
restriction of the determinant to $\MM_S$. This restriction is a polynomial 
in $\R^{|S|}$ whose zero set, if different from $\R^{|S|}$, has dimension 
smaller than $|S|$.
\eproof

\proofof{Theorem~\ref{s2}} 
Case (i): $\MM_S\subseteq\Sigma$. In this case, the desired inequality 
is trivial by Lemma~\ref{s23}. 
\smallskip

Case (ii): $\MM_S\not\subseteq\Sigma$. 
By Lemma~\ref{s22}, with probability $1$, $A$ is non-singular. So,
by Lemma~\ref{s3},
\begin{eqnarray}\label{eq15}
\Prob\{\cd(A)\geq t\}
&=&\Prob\left\{\sum_{i,j\in[n]} |a_{ij}\g_{ji}|\geq t\right\}\notag\\
&=&\Prob\left\{ \sum_{(i,j)\in S}
 \left|\frac{a_{ij}\det(A_{(ij)})}{\det(A)}\right|\geq t\right\}.
\end{eqnarray}
Assume $(1,1)\in S$ and let $x=a_{1}$ be the first column of $A$. 
Also, let $I=\{i\in[n]\mid (i,1)\in S\}$ and  
$x_I$ be the vector obtained by removing
entries $x_i$ with $i\not\in I$. Then,
\begin{equation}\label{eq17}
     x_S\sim N(0,\Id_{|I|}).
\end{equation}
Let $q=\left(\det A_{(11)}, \ldots,\det A_{(n1)}\right)\rt $ and 
$q_{I}$ be the vector obtained by removing entries $q_i$ with 
$i\not\in I$. Clearly, $q_I$ is independent of $x_I$. 
Using this notation, the expansion by the first column yields 
$$
   \det(A)=\sum_{i\in [n]}(-1)^{i+1}a_{i1}\det(A_{(i1)})=x_I\rt q_I.
$$
In addition, $a_{11}\det(A_{(11)})=x_I\rt (q_1e_1)$ 
where $e_1$ is the vector with the first entry equal to $1$ and all
others equal to $0$. Hence, 
$$
     \frac{a_{11}\det(A_{(11)})}{\det(A)}
    =\frac{x_I\rt (q_1e_1)}{x_I\rt q_I}
$$
Using (\ref{eq17}) and Lemma~\ref{sTail} 
(with $p=(q_1e_1)$ and $q=q_I$) we obtain, for $z\geq 2$, 
$$
  \Prob\left\{\left|\frac{a_{11}\det(A_{(11)})}{\det(A)}\right|
  \geq z\right\} \leq \frac{1}{z}. 
$$
The same bound can be proven for all $(i,j)\in S$. 
Using these bounds with $z=\frac{t}{|S|}$ and (\ref{eq15}) we obtain,
\begin{equation}\tag*{\qed}
\Prob\{\cd(A)\geq t\}
 \leq \sum_{(i,j)\in S} \Prob\left\{ 
 \left|\frac{a_{ij}\det(A_{(ij)})}{\det(A)}\right|\geq \frac{t}{|S|}\right\}
 \leq \frac{|S|^2}{t}. 
\end{equation}

The proof of Corollary~\ref{cor:s2} follows from the following 
result by taking $Z=\cd(A)$ and $t_0=|S|^2$. 

\begin{proposition}\label{prop:E}
Let $t_0>0$ and $Z\geq 1$ be a random variable satisfying 
that $\Prob\{Z\geq t\}\leq t_0 t^{-1}$ for all $t\ge t_0$.
Then $\bfE(\log Z)\leq \log t_0 +\frac{1}{\ln b}$ where $b>1$ is the 
base of the logarithm.
\end{proposition}

\proof
We have $\Prob\{\log Z\geq t\}\leq t_0 b^{-t} = b^{-(t-\log t_0)}$
for all $t>\log t_0$. Therefore, 
\begin{equation}\tag*{\qed}
 \bfE(\log Z) = \int_0^{\infty}\Prob\{\log Z\geq s\} ds 
  \leq  \log t_0 + \int_{\log t_0}^{\infty} b^{-(t-\log t_0)} dt
  = \log t_0 + \frac{1}{\ln b}.
\end{equation}

\section{Matrix inversion}

We now focus on the problem of inverting a matrix $A$ and its 
componentwise condition number $\c^{\dagger}(A)$. Our main results in this 
section are the following. 

\begin{theorem}\label{th3}
Let $S\subset [n]^2$ be such that $\MM_S\not\subseteq\Sigma$. Then, 
for all $t\geq 2|S|$, 
$$
  \Prob\{\c^{\dagger}(A)\geq t\}\leq 
  4|S|^2 n^2\frac{1}{t}
$$
where $\Prob$ denotes probability over $A\in\RR_S$. 
\end{theorem}

\begin{corollary}\label{cor:th3}
Let $S\subset [n]^2$ be such that $\MM_S\not\subseteq\Sigma$. Then, 
$$
  \bfE(\log_+(\c^{\dagger}(A)))\leq 2\log n + 2\log |S| 
  +\log 4 +\frac{1}{\ln b}
$$
where $\bfE$ denotes expectation over $A\in\RR_S$. 
\eproof
\end{corollary}

\begin{remark}
Note that, for all monotonic norm on $\MM_S$, 
the bound above also holds for $\m(A)$ by Proposition~\ref{prop:monot}. 
This is in contrast with the lower bound linear in $n$ for the 
expected value of the logarithm of normwise condition numbers, etc. 
\end{remark}

Definition~\ref{def:CN} yields expressions for the (mixed and componentwise) 
condition numbers of a matrix $A$ by taking $\DD=\MM\setminus\Sigma$ and 
$F:\MM\setminus\Sigma\to\MM$ given by $F(A)=A^{-1}$.
For $k,\ell\in[n]$ such that $\g_{k\ell}\neq 0$, we let
$$
 \ckl(A) = \lim_{\delta\rightarrow 0} \sup_{A'\in \BB(A, \delta)}
         \frac{|\g'_{k\ell}-\g_{k\ell}|}{\delta|\g_{k\ell}|}
$$
and for $k,\ell\in[n]$ such that $\g_{k\ell}= 0$, we let $\ckl(A) = 0$ if
$$
\lim_{\delta\rightarrow 0} \sup_{A'\in \BB(A, \delta)}
     \frac{|\g'_{k\ell}-\g_{k\ell}|}{\delta}=0
$$
and $\ckl(A)=\infty$ otherwise. Then 
$$
   \c^{\dagger}(A)=\max_{k,\ell\in[n]}\ckl(A).
$$
Similarly, for a norm $\|\ \|$ on $\MM$, 
$$
 \m(A)= \lim_{\delta\rightarrow 0} \sup_{A'\in \BB(A, \delta)}
      \frac{\|\Gamma'-\Gamma\|}{\delta\|\Gamma\|}.
$$

\begin{lemma}\label{s1}
For $A\in\MM\setminus\Sigma$ and 
$k,\ell\in[n]$,
$$
   \ckl(A) \leq \cd(A) + \cd(A_{(\ell k)}).
$$
\end{lemma}

\proof
We divide the proof by cases. Case (i): $\g_{k\ell}\neq 0$. 

Let $\delta>0$ be sufficiently small so that 
$\BB(A,\delta)\cap\Sigma=\emptyset$ and, for all $A'\in\BB(A,\delta)$,  
$\left|\frac{\det(A')-\det(A)}{\det(A)}\right|<1$. 
Let $A'\in\BB(A,\delta)$. 

Since $\g_{k\ell}=\frac{\det(A_{(\ell k)})}{\det(A)}$, 
\begin{eqnarray*}
\frac{\g'_{k\ell}-\g_{k\ell}}{\g_{k\ell}}
&=&\frac{\det(A)}{\det(A_{(\ell k)})}\left(\frac{\det(A'_{(\ell k)})}{\det(A')}-\frac{\det(A_{(\ell k)})}{\det(A)}\right)\\
&=&\frac{\det(A)}{\det(A_{(\ell k)})}\frac{\det(A'_{(\ell k)})}{\det(A')}-1\\
&=&\frac{1+\frac{\det(A'_{(\ell k)})-\det(A_{(\ell k)})}{\det(A_{(\ell k)})}}{1+\frac{\det(A')-\det(A)}{\det(A)}}-1\\
&=&\frac{\frac{\det(A'_{(\ell k)})-\det(A_{(\ell k)})}{\det(A_{(\ell k)})}-\frac{\det(A')-\det(A)}{\det(A)}}{1+\frac{\det(A')-\det(A)}{\det(A)}}.\\
\end{eqnarray*}
Using that $\left|\frac{\det(A')-\det(A)}{\det(A)}\right|<1$, 
\begin{eqnarray*}
  \left|\frac{\g'_{k\ell}-\g_{k\ell}}{\g_{k\ell}}\right|
 &\leq& \frac{\left|\frac{\det(A'_{(lk)})-\det(A_{(\ell k)})}
        {\det(A_{(\ell k)})}\right|
       +\left|\frac{\det(A')-\det(A)}{\det(A)}\right|}
        {1-\left|\frac{\det(A')-\det(A)}{\det(A)}\right|}
\end{eqnarray*}
and therefore
\begin{eqnarray*}
  \sup_{A'\in \BB(A, \delta)}& &
  \kern-25pt
  \left|\frac{\g'_{k\ell}-\g_{k\ell}}{\delta \g_{k\ell}}\right|\\
 &\leq & 
  \frac{\sup_{A'\in\BB(A, \delta)}
  \left|\frac{\det(A'_{(\ell k)})-\det(A_{(\ell k)})}
  {\delta\det(A_{(\ell k)})}\right|+\sup_{A'\in \BB(A, \delta)}
  \left|\frac{\det(A')-\det(A)}{\delta\det(A)}\right|}
  {1-\sup_{A'\in \BB(A, \delta)}
  \left|\frac{\det(A')-\det(A)}{\det(A)}\right|}.
\end{eqnarray*}
Taking limits for $\delta \rightarrow 0$ on both sides we get
$$
\ckl(A)\leq \cd(A) + \cd(A_{(\ell k)}).
$$
Case (ii): $\g_{k\ell}=0$ and
$$
\lim_{\delta\rightarrow 0} \sup_{A'\in \BB(A, \delta)}
    \frac{|\g'_{k\ell}|}{\delta}=0.
$$
In this case, $\ckl(A)=0$ and the statement holds. 

Case (iii): $\g_{k\ell}= 0$ and 
$$
0 \neq \lim_{\delta\rightarrow 0} \sup_{A'\in \BB(A, \delta)}
   \frac{|\g'_{k\ell}|}{\delta}
 =\lim_{\delta\rightarrow 0} \sup_{A'\in \BB(A, \delta)}
 \frac{|\det(A'_{\ell k})|}{\delta|\det(A')|}.
$$
In this case 
$$
  \lim_{\delta\rightarrow 0} \sup_{A'\in \BB(A, \delta)}
  \frac{|\det(A'_{\ell k})|}{\delta}\neq 0
$$
and therefore $\cd(A'_{\ell k})=\infty$. The statement 
holds as well.
\eproof

\proofof{Theorem~\ref{th3}}
By definition of $\c^{\dagger}(A)$,
$$
\Prob\{\c^{\dagger}(A)\geq t\}
=\Prob\left\{\max_{k,\ell\in[n]}\ckl(A)\geq t\right\}
\leq \sum_{k,\ell\in[n]}\Prob\{\ckl(A)\geq t\}.
$$
By Lemma~\ref{s22}, with probability $1$, $A$ is non-singular. So,
we can apply Lemma~\ref{s1} to obtain
\begin{eqnarray*}
\Prob\{\ckl(A)\geq t\} 
&\leq&\Prob\left\{\cd(A)\geq \frac{t}{2}\right\} +
   \Prob\left\{\cd(A_{(k\ell)})\geq \frac{t}{2}\right\} \\
&\leq& 4|S|^2\frac{1}{t}
\end{eqnarray*}
the last inequality by Theorem~\ref{s2}. The statement now follows.
\eproof

\section{Linear Equations Solving}

We finally deal with the problem of solving linear systems of equations. 

\begin{theorem}\label{th5}
Let $S\subset [n]^2$ be such that 
$\MM_S\not\subseteq\Sigma$. Then, for all $t\geq 2(|S|+n)$, 
$$
  \Prob\{\c(A,b)\geq t\}\leq 10|S|^2 n\frac{1}{t} 
$$
where $\Prob$ denotes probability over $(A,b)\in\RR_S\times N(0,\Id_n)$. 
\end{theorem}

\begin{corollary}\label{cor:th5}
Let $S\subset [n]^2$ be such that 
$\MM_S\not\subseteq\Sigma$. Then,
\begin{equation}\tag*{\qed}
  \bfE(\log_+(\c(A,b)))\leq \log n + 2\log |S| + 
   \log 10 +\frac{1}{\ln b}.
\end{equation}
\end{corollary}

Definition~\ref{def:CN} yields again expressions for the (mixed and
componentwise) condition numbers of a pair $(A,b)$ by taking
$\DD=(\MM\setminus\Sigma)\times\R^n$ and
$F:(\MM\setminus\Sigma)\times\R^n\to\R^n$ given by $F(A,b)=A^{-1}b$.

For $A\in\MM\setminus\Sigma$ and $b\in\R^n$ we denote $x=A^{-1}b$. 
For $k\in[n]$ such that $x_k\neq 0$ we let
$$
  \ck(A,b) = \lim_{\delta\rightarrow 0}\sup_{(A',b')\in \BB((A,b), \delta)}
                \frac{|x'_k-x_k|}{\delta|x_k|}.
$$
For $k\in[n]$ such that $x_k= 0$ we let $\ck(A, b) = 0$ if
$$
  \lim_{\delta\rightarrow 0} \sup_{(A',b')\in \BB((A,b),\delta)}
                    \frac{|x'_k-x_k|}{\delta}=0
$$
and $\ck(A,b)=\infty$ otherwise. Then
$$
  \c(A,b) = \max_{k\in[n]}\ck(A,b).
$$
Similarly, for a norm $\|\ \|$ in $\R^n$, 
$$
  \m(A,b) = \lim_{\delta\rightarrow 0} \sup_{(A',b')\in \BB((A,b),\delta)}
    \frac{\|x'-x\|}{\delta\|x\|}.
$$

In what follows let $R_k$ be the matrix obtained by replacing the $k$th column 
of $A$ by $b$.

\begin{lemma}\label{th6}
For any non-singular matrix $A$ and $k\in[n]$,
$$
 \ck(A,b)\leq \cd(A)+\cd(R_k).
$$
\end{lemma}

\proof
By Crammer's rule,
$$
   x_k=\frac{\det(R_k)}{\det(A)}.
$$
The rest of this proof is similar to the proof of Lemma~\ref{s1}.
\eproof

\proofof{Theorem~\ref{th5}}
It follows the lines of that of Theorem~\ref{th3}. First, we get 
$$
  \Prob\{\c(A,b)\geq t\}\leq \sum_{k\in[n]}\Prob\{\ck(A,b)\geq t\}.
$$
Then, we apply Lemma~\ref{th6} (using that, with probability $1$, 
$A\not\in\Sigma$) and the fact that $|S|\geq n$ to get
\begin{eqnarray*}
 \Prob\{\ck(A,b)\geq t\} &\leq& 
 \Prob\left\{\cd(A)\geq \frac{t}{2}\right\} 
 +\Prob\left\{\cd(R_k)\geq \frac{t}{2}\right\}\\ 
 &\leq& 2|S|^2 \frac{1}{t} 
       + 2(|S|+n)^2\frac{1}{t}
  \leq  10 |S|^2 \frac{1}{t}
\end{eqnarray*}
from which the statement follows.
\eproof

\section{Additional Remarks}\label{sec:final}

\noindent
{\bf (1)\quad}
To obtain bounds for the average loss of precision of random 
triangular systems one may combine Theorem~\ref{th5} with the following 
result by Wilkinson~\cite[Ch.3,\S19]{Wilkinson63} which we quote as given 
in~\cite{Higham:89}.

\begin{theorem}
Let $T\in\R^{n\times n}$ be a nonsingular triangular matrix, an assume 
$nu<0.1$ (here $u$ is the round-off unit). Then, the computed solution 
$\hat x$ to the system $Tx=b$ satisfies 
$$
   (T+E)\hat x=b,
$$
where, for some universal constant $c$,  
\begin{equation}\tag*{\qed}
   |e_{ij}|\leq (|i-j|+2)cu|t_{ij}|. 
\end{equation}
\end{theorem}
The use of $\c(A)$ actually yields average loss of precision with the 
latter measured componentwise in the computed solution. 
\medskip

\noindent
{\bf (2)\quad}
The bound in Corollary~\ref{cor:th3} appears to be worse than what computer 
simulations suggest $\bfE(\log \c^{\dagger}(L_n))$ should be. 
In~\cite{CuWe:07} matrices $L_n$ were generated for various values of 
$n$ and an experimental mean of $\bfE(\log \c^{\dagger}(L_n))$ was 
obtained from these values. A linear regression for these means shows 
a best fit of $3.065 \log n -1.1466$. A probable source of (a good part 
of) the difference of this value with the estimate  
$6\log n + \Oh(1)$ 
following from Corollary~\ref{cor:th3} 
is the broad bound 
$\Prob\left\{\max_{k,\ell\in[n]}\ckl(A)\geq t\right\}
\leq \sum_{k,\ell\in[n]}\Prob\{\ckl(A)\geq t\}$
in the proof of Theorem~\ref{th3}. In addition to this, the bound 
$\bfE(\m^{\dagger}(L_n))\leq  \bfE(\c^{\dagger}(L_n))$ following 
from Proposition~\ref{prop:monot} may be coarse as well. Numerical 
experiments in~\cite{CuWe:07} suggest a best fit of 
$\bfE(\log \m^{\dagger}(L_n))\approx 1.5334 \log n -0.5723$. 
\medskip

\noindent
{\bf (3)\quad}
Section~6 of~\cite{ViTr98} discusses stability of Gaussian elimination. 
Having shown that, almost surely, $\kappa(L_n)\approx 2^n$ the authors 
reflect on how this behavior can be reconciled with the fact that 
``Gaussian elimination is overwhelmingly stable.'' They point out that 
``The reason appears to be statistical: the matrices $A$ for which 
$\|L^{-1}\|$ is large occupy an exponentially small proportion of the 
space of matrices'' a claim for which experimental evidence is given 
in~\cite{TreSch90}. It would thus follow that ``{\em the matrices $L$ produced by 
Gaussian elimination are far from random.}'' 

Our results show that, in addition, Gaussian elimination needs much less 
than producing matrices $L$ in the vanishingly small set of triangular 
matrices with $\kappa(L)$ small. It is enough to produce matrices $L$ 
{\em outside} the vanishingly small set of matrices with $\c(L,b)$ large. 
\medskip

\noindent
{\bf Aknowledgement.}\quad We are grateful to Ernesto Mordecki for helpful 
discussions. 

{\small

}

\end{document}